\documentclass{amsart}
\begin{document}
\title[The finite representation property fails]%
{The finite representation property fails for composition and intersection}
\author{Roger D.\ Maddux}
\address{Department of Mathematics
\\Iowa State University 
\\Ames, Iowa 50011-2066, USA}
\email{ maddux@iastate.edu }
\date{Written May 14, 2013; revised June 25, 2014}
\begin{abstract}
The title theorem is proved by example: an algebra of binary relations, closed under
intersection and composition, that is not isomorphic to any such algebra on a finite
set.
\end{abstract}
\maketitle
\def\rp{\mathop{;}}
\def\A{\mathcal{A}}
\def\l{r}
\def\0{z}
\def\e{e}
\def\Q{\mathbb{Q}}
\def\lo{\l\cdot\overline\0}
\newtheorem{theorem}{{\bf Theorem}}
Let $K$ be a class of algebras for which there is a notion of ``representability over a
set $U$''. That is, for every set $U$, some algebras of $K$ are said to be {\bf
representable over $U$}, while others are not. We say that $K$ has the {\bf finite
representation property} if every finite algebra in $K$ that has a representation over
\emph{some} set has a representation over a \emph{finite} set.

$K$ may be defined abstractly, as a class of algebras of some particular similarity
type, satisfying some conditions which, if they are all universally quantified
equations, means that $K$ is a variety. In this case some definition of representability
is still required. However, if $K$ is taken to be a class of algebras described in some
concrete set-theoretical manner, then we may wish representability to simply be
membership in $K$. An example of this type, one that fails to have the finite
representation property, is considered here.\footnote{Theorem~\ref{th2} was concocted
during conversations with Jeremy Alm, Robin Hirsch, Richard Kramer, and Andy Ylvisaker,
at Iowa State University in early May of 2013.}

Let $K$ be the class of algebras of the form $(A,\rp,\cdot)$, where $\rp$ and $\cdot$
are binary operations on $A$, such that, for some set $U$, $A$ is a set of binary
relations on $U$, and for all $a,b\in A$, $a\rp b$ is the compositum of the relations
$a$ and $b$, in that order, while $a\cdot b$ is the intersection of $a$ and $b$ (in
either order). In more detail, for all $a,b\in A$ we have
\begin{align*}
	a\rp b	&=\{(x,y):\text{for some }z\in U,\,(x,z)\in a\text{ and }(z,y)\in b\},
\\	a\cdot b&=\{(x,y):(x,y)\in a\text{ and }(x,y)\in b\}.
\end{align*}
An algebra in $K$ can be described simply as a set of relations (on some base set $U$)
that is closed under composition and intersection. Every algebra in $K$ is representable
over some set, namely, the base set $U$ used to specify the algebra, which may be
necessarily infinite. 
\begin{theorem}\label{th1}
	$K$ does \emph{not} have the finite representation property.
\end{theorem}
We will show this by giving an example of an algebra $\A$ in $K$ that is not isomorphic
to any algebra in $K$ with a finite base set. The example is called the {\bf point
algebra} (by analogy with the relation algebra having the same name). The base set of
$\A$ is the set $\Q$ of rational numbers, and the elements of $\A$ are these three
relations:
\begin{align*}
	\l&:=\{(x,y):x,y\in\Q\land x<y\}
\\	\0&:=\emptyset
\\	\e&:=\{(x,x):x\in\Q\}
\end{align*}
The tables for the two operations are given below, with the entries that are actually
used later enclosed in boxes:
$$\def\bl{\boxed{\l}}
\def\bo{\boxed{\0}}
\begin{array}{|c|ccc|}		\hline
	\rp	&\0	&\e	&\l	\\\hline
	\0	&\0	&\0	&\bo	\\
	\e	&\0	&\e	&\bl	\\
	\l	&\bo	&\bl	&\bl	\\\hline
\end{array}\qquad
\begin{array}{|c|ccc|}		\hline
	\cdot	&\0	&\e	&\l	\\\hline
	\0	&\0	&\0	&\0	\\
	\e	&\0	&\e	&\bo	\\
	\l	&\0	&\bo	&\l	\\\hline
\end{array}
$$
The structure of $\A$ is completely specified by the two tables, and the second table is
determined entirely by either of its boxed entries. That $\A$ belongs to $K$ follows
from the fact that if the elements $\l,\e,\0$ are defined as the binary relations given
above, then the two tables can be deduced from the definitions. What we do next is
assume that $\A$ has a representation over some set $U$, and show that $U$ must be
infinite.
\begin{theorem}\label{th2}
If $U$ is a non-empty set with distinct relations $\0,\l,\e\subseteq U\times U$
satisfying $\l\rp\e=\l=\e\rp\l$, $\l\rp\l=\l$, $\0\rp\l=\0=\l\rp\0$, and $\l\cdot\e=\0$,
then $U$ is infinite.
\end{theorem}
\proof
First we show the intersection of the identity relation on $U$ with $\l$ is included in
$\0$, that is,
\begin{equation}\label{6}
	Id_U\cdot\l\subseteq\0.
\end{equation}
To show this, we assume
\begin{equation}\label{1}
	(x,x)\in\l
\end{equation}
and derive $(x,x)\in\0$. From \eqref{1} and $\l=\l\rp\e$ we get
$(x,x)\in\l\rp\e$, hence we know there is some $y\in U$ such that
\begin{equation}\label{2}
	(x,y)\in\l,
\end{equation}
\begin{equation}\label{3}
	(y,x)\in\e.
\end{equation}
From \eqref{3} and \eqref{1} we get $(y,x)\in\e\rp\l$, but $\e\rp\l=\l$, so
\begin{equation}\label{4}
	(y,x)\in\l.
\end{equation}
Then \eqref{4} and \eqref{3} give us $(y,x)\in\l\cdot\e$, but $\l\cdot\e=\0$, so
\begin{equation}\label{5}
	(y,x)\in\0.
\end{equation}
From \eqref{2} and \eqref{5} we have $(x,x)\in\l\rp\0$, but $\l\rp\0=\0$, so
$(x,x)\in\0$. This completes the proof of \eqref{6}.

Note that~\eqref{6} is equivalent to $\l\cdot\overline\0\subseteq\overline{Id_U}$,
{\it~i.e.}, the intersection of $\l$ with the complement of $\0$ (with respect to
$U\times U$) is a diversity relation (included in the complement of the identity
relation on $U$).  Note also that $\0\subseteq\l$ and $\0\subseteq\e$ because
$\l\cdot\e=\0$.  All three relations $\0,\e,\l$ must be distinct, for otherwise we do
not have a representation, hence $\lo\neq\emptyset\neq\e\cdot\overline\0$.  Since $\lo$
is a non-empty diversity relation, there are distinct $x_0,y\in U$ such that
\begin{align}
\label{7}	(x_0,y)\in\l,	\\
\label{8}	(x_0,y)\in\overline\0.
\end{align}
From \eqref{7} and $\l=\l\rp\l$ we know there is some $x_1\in U$ such that
\begin{align}
\label{10}	(x_0,x_1)\in\l,	\\
\label{11}	(x_1,y)\in\l.
\end{align}
If $(x_0,x_1)\in\0$ then $(x_0,y)\in\0\rp\l$ by \eqref{11}, but $\0\rp\l=\0$, so we get
$(x_0,y)\in\0$, contradicting~\eqref{8}. Therefore $(x_0,x_1)\in\overline\0$, hence
$(x_0,x_1)\in\lo$ by \eqref{7}.  Similarly, if $(x_1,y)\in\0$ then $(x_0,y)\in\l\rp\0$
by \eqref{10}, but $\l\rp\0=\0$, so we get $(x_0,y)\in\0$, contradicting~\eqref{8}.
Therefore $(x_1,y)\in\overline\0$.  

So far we have in fact proved that $\lo$ is a non-empty dense diversity relation: there
are distinct $x_0,x_1,y\in U$ such that $(x_0,y),(x_0,x_1),(x_1,y)\in\lo$. We have also
achieved the first stage (with $n=1$) in the construction of $y$, $x_0$, $x_1$,
$x_2$,~\dots, $x_n$ such that
\begin{equation}\label{14}
	(x_i,x_j)\in\lo	\quad\text{whenever $0\leq i<j\leq n$,}
\end{equation}
\begin{equation}\label{15}
	(x_i,y)\in\lo	\quad\text{whenever $0\leq i\leq n$.}
\end{equation}
We continue this construction through one more stage. Apply the density of $\lo$ to the
assumption $(x_n,y)\in\lo$, obtaining some $x_{n+1}$ such that
\begin{equation}\label{a}
	(x_n,x_{n+1})\in\lo,
\end{equation}
\begin{equation}\label{b}
	(x_{n+1},y)\in\lo.
\end{equation}
Obviously \eqref{b} implies that \eqref{15} holds with $n+1$ in place of $n$.  To see
the same for \eqref{14}, let $0\leq i<j\leq n+1$.  If $j<n+1$ we are done,
by~\eqref{14}, so we may assume $j=n+1$.  We wish to show $(x_i,x_{n+1})\in\lo$.  This
holds by \eqref{a} if $i=n$, so assume $i<n$.  We have $(x_i,x_n)\in\l$ by \eqref{14}
and $(x_n,x_{n+1})\in\l$ by \eqref{a}, so $(x_i,x_{n+1})\in\l\rp\l=\l$.  If
$(x_i,x_{n+1})\in\0$ then $(x_i,y)\in\0\rp\l$ by \eqref{b}, but $\0\rp\l=\0$, so
$(x_i,y)\in\0$, contradicting~\eqref{15}, hence $(x_i,x_{n+1})\in\overline\0$. Thus we
have $(x_i,x_{n+1})\in\lo$. This construction may be continued indefinitely, so $U$ must
be infinite.
\endproof
\end{document}